\documentclass[a4paper,10pt]{article}

\usepackage{amsmath, amsfonts, amssymb, latexsym, mathrsfs, graphicx}

\newcommand{\mmap}[3]{\ar@/#2/[#1]|*=0{\rotatebox{#3}{$\scriptsize |$}}} 
\newcommand{\picar}[3]{\ar@{}[#1]|*{\rotatebox{#2}{$\scriptsize #3$}}} 

\newcommand{\sa}{{\mathscr{A}}}

\newcommand{\scra}{{\mathscr{P}(\mathscr{A})}}
\newcommand{\sv}{{\mathscr{V}}}
\newcommand{\ssc}{{\mathscr{C}}}
\newcommand{\oop}{\operatorname{op}}
\newcommand{\op}{^{\oop}}

\usepackage[all]{xypic}

\begin{document}

\begin{center}
{\underline{{\Large On The Existence Of Category Bicompletions}}}
\newline

{{\small Brian J. DAY}}

{{\small May 4, 2009.}}
\end{center}

\noindent \emph{\small \underline{Abstract:} A completeness conjecture is advanced concerning the free small-colimit completion $\scra$ of a (possibly large) category $\sa$.
The conjecture is based on the existence of a small generating-cogenerating set of objects in $\sa$.
We sketch how the validity of the result would lead to the existence of an Isbell-Lambek bicompletion $\mathscr{C}(\sa)$ of such an $\sa$, without a ``change-of-universe'' procedure being necessary to describe or discuss the bicompletion}
\newline

\begin{center}
\begin{tabular}{c}\hline
\phantom{a}\phantom{a}\phantom{a}\phantom{a}\phantom{a}\phantom{a}\phantom{a}\phantom{a}\phantom{a}\phantom{a}\phantom{a}\phantom{a}\phantom{a}\phantom{a}\phantom{a}\phantom{a}
\end{tabular}
\end{center}

All categories, functors, and natural transformations, etc., shall be relative to a basic complete and cocomplete symmetric monoidal closed category $\sv$ with all intersections of subobjects.
A tentative conjecture, based partly on the results of [3], is that if $\sa$ is a (large) category containing a small generating and cogenerating set of objects, then $\scra$ (which is the free small-colimit completion of $\sa$ with respect to $\sv$) is not only cocomplete (as is well known), but also complete with all intersections of subobjects.
\newline

If this conjecture is true, then one can establish the existence of a resulting ``Isbell-Lambek'' bicompletion of such an $\sa$, along the lines of [1] \S 4, using the Yoneda embedding $Y: \sa\subset\mathscr{P}(\sa)$.
This proposed bicompletion, denoted here by $\mathscr{C}(\sa)$, has the same ``size'' as $\sa$ and is, roughly speaking, the (replete) closure in $\mathscr{P}(\sa)$, under both iterated limits and intersections of subobjects, of the class (i.e. large set) of all representable functors from $\sa\op$ to $\sv$.
\newline

More precisely, one can construct $\mathscr{C}(\sa)$ directly using the Isbell-conjugacy adjunction

$$\xymatrix{ \mathscr{P}(\sa) \ar@<+1ex>[rr]^{Lan_{Y}(Z)} && \ar@<+1ex>[ll]^{R} \mathscr{P}(\sa\op)\op }$$

\noindent whose existence (see [3] \S 9) follows from the conjectured completeness of both $\mathscr{P}(\sa)$ and $\mathscr{P}(\sa\op)$, and where $Z:\sa\to\mathscr{P}(\sa\op)\op$ is the dual of the Yoneda embedding $\sa\op\subset\mathscr{P}(\sa\op)$.
Thus we proceed by factoring the left adjoint $Lan_{Y}(Z)$ as a reflection followed by a conservative left adjoint

$$\xymatrix{ \ar@/_/@<-2.5ex>[d]  \mathscr{P}(\sa) \ar@<+0.8ex>[rr]^-{Lan_{Y}(Z)} && \ar@<+0.8ex>[ll]^-{R} \mathscr{P}(\sa\op)\op \ar[lld]  \\   \ssc(\sa) \ar@{}[u]|*=0{\cup}  \ar@<-1.5ex>[rru]_-{cons}  && .  }$$

\noindent Such a factorization exists by [1] Theorem 2.1 and is essentially unique by [1] Proposition 5.1. 
Moreover, the induced full embedding

$$\sa\subset\ssc(\sa)$$

\noindent then preserves any small limit or small colimit that already exists in $\sa$.
\newline

One important consequence is that various results from [2] on monoidal biclosed completion of categories can be accordingly revamped using such a bicompletion $\ssc(\sa)$; see also  [3] \S 7, which describes some examples where $\mathscr{P}(\sa)$ is monoidal or monoidal biclosed.
Note that here especially one could conveniently avoid the awkward ``change-of-$\sv$-universe'' procedure employed in [2].
\newline
\newline

\begin{center}
{\underline{{\Large References.}}}
\end{center}

\begin{center}
\begin{tabular}{l l}
\textrm{[1]} & B. J. Day, ``On Adjoint-Functor Factorization", \\
 & Lecture Notes in Mathematics, 420 (Springer-Verlag 1974), Pg. 1-19.\\
 & \\
\textrm{[2]} & B. J. Day, ``On Closed Categories Of Functors II'', \\
 & Lecture Notes in Mathematics, 420 (Springer-Verlag 1974), Pg. 20-54.  \\
 & \\
\textrm{[3]} & B. J. Day and S. Lack, ``Limits Of Small Functors'' \\ 
 & J. Pure Appl. Alg., 210 (2007), Pg. 651-663. \\
 & \\
 & \\
 & \\
 & \\
\end{tabular}

\small{Mathematics Dept., Faculty of Science, Macquarie University, NSW 2109, Australia.}
\newline

Any replies are welcome through Tom Booker (thomas.booker@students.mq.edu.au), who kindly typed the manuscript.

\end{center}

\end{document}